\documentclass[11pt]{article}
\usepackage{xcolor}
\usepackage{amsmath}
\usepackage{amssymb}
\usepackage{latexsym}
\usepackage[all,cmtip]{xy}
\usepackage{graphicx, float, color}
\usepackage{accents}

\newtheorem{thm}{Theorem}
\newtheorem{definition}[thm]{Definition}
\newtheorem{lemma}[thm]{Lemma}

\begin{document}

\title{Moves on Filtered PL Manifolds and Stratified PL Spaces}

\author{Louis Crane and David N.\ Yetter}

\maketitle

\begin{abstract}
We extend results of Pachner \cite{Pa87}, \cite{Pa91} and Casali \cite{Ca95} to give finite sets of moves relating triangulations of PL manifolds respecting filtrations by locally flat manifolds and stratifications in which a finite family of simple local models exists for neighborhoods of strata.
\end{abstract}

Recent work of F\"{u}chs et al. \cite{FSV13a}, \cite{FSV13b} on topological quantum field theories with defects prompts the consideration of adequate combinatorial moves to describe simplicial constructions after the manner of Turaev and Viro \cite{TV92} in 3-dimensions or Crane and Yetter \cite{CY93} in 4-dimensions in a setting in which manifolds are equipped with a stratification by ``defects'' in various codimensions.  Plainly Alexander moves will suffice, but as in the defect-free case, these are difficult to use because there are infinitely many combinatorial types once the dimension exceeds two.

The goal of this work is to establish finite families of moves after the manner of Pachner \cite{Pa87}, \cite{Pa91}, which are adequate to relate triangulations of PL manifolds filtered by locally flat submanifolds or with stratifications for which a finite family of simple local models exists for neighborhoods of strata.

Throughout we assume a familiarity with the basic notions of piecewise linear topology, including triangulations, combinatorial manifolds, various types of subdivisions, stars, links, and joins.  The reader needing an introduction to or refresher on these matters is referred to the classic book by Rourke and Sanderson \cite{RS72}.  We also assume all manifolds and manifolds with boundary are compact, and that triangulations of a space are those for a fixed PL structure.

We begin by recalling the main theorems of the work on which we build:

\begin{thm} [Pachner \cite{Pa87}, \cite{Pa91} (also \cite{Li99} Theorem 5.9)]  Two closed combinatorial $n$-man\-i\-folds are piecewise linear homeomorphic if and only if they are bistellar equivalent.
\end{thm}

\noindent and

\begin{thm} [Casali \cite{Ca95} MAIN THEOREM] \label{Casali} Let $M$ and $M^\prime$ be two PL-manifolds, and let $K$ (resp. $K^\prime$) be a simplicial triangulation of $M$ (resp. $M^\prime$) with $\partial K = \partial K^\prime$.  Then $M$ and $M^\prime$ are PL-homeomorphic if and only if $K$ and $K^\prime$ are bistellar equivalent.
\end{thm}

Recall the definition of bistellar equivalence:  

\begin{definition} [Pachner \cite{Pa87}]  Let $A$ be a $k$-simplex ($0\leq k \leq n$) of ${\rm Int}K = K \setminus \partial K$ such that its link ${\rm lk}(A; K)$ is $\partial B$ for an $(n-k)$-simplex $B$ which is not in $K$.  Then the {\em bistellar $k$-operation} $\chi_{(A,B)}$ on $K$ is the process of replacing the star of $A$, ${\rm st}(A; K) = A\ast \partial B$ with $\partial A \ast B$, that is

\[ \chi_{(A,B)}K = (K \setminus A \ast \partial B) \cup (\partial A \ast B) . \]

\noindent Two triangulations of a PL-manifold $M$ are {\em bistellar equivalent} if one is obtained from the other by a finite sequence of bistellar $k$-operations (for various $k$).
\end{definition}

The key lemma in Casali's proof of \ref{Casali} is 

\begin{lemma} [Casali \cite{Ca95} Proposition 4] If $A \in {\rm Int}K$, then $K$ and the stellar subdivision of $K$ at (any point $a$ in) $A$, $\sigma_A K$ are bistellar equivalent.
\end{lemma}

A careful reading of Casali's proof shows that, in fact, it carries a stronger result which is the basis for the results herein:

\begin{lemma} \label{BetterC'sLemma} If $K$ is a triangulation of a manifold with boundary, $M$, and $L$ is any subcomplex of $K$ which contains $\partial M$, and $A$ is any simplex of $K \setminus L$, then $K$ and $\sigma_A K$ are bistellar equivalent by bistellar moves which leave $L$ unchanged.
\end{lemma}

\noindent{\sc Proof.} The proof is Casali's, with the condition of internalness replaced throughout with the condition of not lying in $L$. $\Box$

Our first object of interest is given by

\begin{definition} A {\em locally flat filtration by submanifolds} of an $n$-manifold $M$ is a sequence of (possibly empty) submanifolds

\[ M_0 \subset M_1 \subset \ldots \subset M_{n-1} \subset M_n = M \]

\noindent in which $M_k$ is $k$-dimensional, and each inclusion $\iota_{k,\ell}:M_k \hookrightarrow M_\ell$ for $k < \ell$ is locally flat.  For brevity, we will call a manifold equipped with a locally flat filtration by submanifolds a {\em filtered manifold}.  By a {\em triangulation of a filtered manifold $M$} we mean a triangulation of $M$ which restricts to a triangulation of each $M_k$.
\end{definition}

Note that a filtered manifold is a very simple sort of stratified space, with the codimension $k$ stratum given by $M_{n-k} \setminus M_{n-k-1}$.

Now, it is well known that local flatness of a submanifold implies the existence of a collar neighborhood (cf. \cite{RS72}). So for every filtered manifold we have a diagram of PL inclusions
\bigskip

\hspace{2.8cm}
\setlength{\unitlength}{3947sp}%
\begingroup\makeatletter\ifx\SetFigFont\undefined%
\gdef\SetFigFont#1#2#3#4#5{%
  \reset@font\fontsize{#1}{#2pt}%
  \fontfamily{#3}\fontseries{#4}\fontshape{#5}%
  \selectfont}%
\fi\endgroup%
\begin{picture}(3105,615)(1186,-751)
\put(1201,-286){\makebox(0,0)[lb]{\smash{{\SetFigFont{12}{14.4}{\rmdefault}{\mddefault}{\updefault}{\color[rgb]{0,0,0}
$M_0 \times I^n \subset M_1 \times I^{n-1} \subset . . . \subset M_{n-1} \times I \subset M_n = M$}%
}}}}
\put(1501,-511){\makebox(0,0)[lb]{\smash{{\SetFigFont{12}{14.4}{\rmdefault}{\mddefault}{\updefault}{\color[rgb]{0,0,0}$\cup$}%
}}}}
\put(2451,-511){\makebox(0,0)[lb]{\smash{{\SetFigFont{12}{14.4}{\rmdefault}{\mddefault}{\updefault}{\color[rgb]{0,0,0}$\cup$}%
}}}}
\put(4276,-511){\makebox(0,0)[lb]{\smash{{\SetFigFont{12}{14.4}{\rmdefault}{\mddefault}{\updefault}{\color[rgb]{0,0,0}$\cup$}%
}}}}
\put(4926,-511){\makebox(0,0)[lb]{\smash{{\SetFigFont{12}{14.4}{\rmdefault}{\mddefault}{\updefault}{\color[rgb]{0,0,0}$\parallel$}%
}}}}
\put(1446,-736){\makebox(0,0)[lb]{\smash{{\SetFigFont{12}{14.4}{\rmdefault}{\mddefault}{\updefault}{\color[rgb]{0,0,0}$M_0 \;\;\; \subset$}%
}}}}
\put(2376,-736){\makebox(0,0)[lb]{\smash{{\SetFigFont{12}{14.4}{\rmdefault}{\mddefault}{\updefault}{\color[rgb]{0,0,0}$M_1 \;\; \subset$ }%
}}}}
\put(3151,-736){\makebox(0,0)[lb]{\smash{{\SetFigFont{12}{14.4}{\rmdefault}{\mddefault}{\updefault}{\color[rgb]{0,0,0}. . .}%
}}}}
\put(3626,-736){\makebox(0,0)[lb]{\smash{{\SetFigFont{12}{14.4}{\rmdefault}{\mddefault}{\updefault}{\color[rgb]{0,0,0}$\subset$}%
}}}}
\put(4116,-736){\makebox(0,0)[lb]{\smash{{\SetFigFont{12}{14.4}{\rmdefault}{\mddefault}{\updefault}{\color[rgb]{0,0,0}$M_{n-1}$ }%
}}}}
\put(4641,-736){\makebox(0,0)[lb]{\smash{{\SetFigFont{12}{14.4}{\rmdefault}{\mddefault}{\updefault}{\color[rgb]{0,0,0}$\subset M_n$}%
}}}}
\end{picture}%

\bigskip

From this it follows from this that if $B_k$ with $\accentset{\circ}{B}_k \subset M_k \setminus M_{k-1}$ is a triangulated ball with interior contained in one of the strata that the inclusion $B_k \hookrightarrow M$ extends to an inclusion $B_k \times I^{n-k} \hookrightarrow M$ with the filtration of $M$ restricting to a filtration of $B_n \times I^{n-k}$ of the form

\[ B_k \subset B_k \times I \subset \ldots \subset B_k \times I^{n-k-1} \subset B_k \times I^{n-k} \;\;\; (\dagger)\]

This local structure will give us the control of the triangulation near a stratum needed to reduce all changes of triangulation to a finite set of moves.

\begin{definition} By a {\em suspension} $\Sigma K$ of a complex $K$ we mean a complex of the form $K \ast \{v_+, v_-\}$ where the set with two elements is understood as a complex whose only simplexes are two distinct vertices.  We denote an m-fold iterated suspension of $K$ by $\Sigma^m K$.  

If $K$ is a subcomplex of $M_k$ in a filtered manifold, a {\em filtered suspension of $K$} is a subcomplex of $M_{k+1}$ which is a supension of $K$ with the vertices $v_+$ and $v_-$ lying in $M_{k+1} \setminus M_k$.  An {\em $m$-fold iterated filtered suspension of $K$ is a subcomplex of the $M_{k+m}$} in which each $\Sigma^{\ell +1}K$ is a filtered supsension of $\Sigma^\ell K$.

An {\em extended bistellar operation} is the operation which replaces $\Sigma^{n-k}(A \ast \partial B)$ with $\Sigma^{n-k}(\partial A \ast B)$, where $A$ is an $\ell$-simplex of the triangulation of $M_k$, $B$ is a simplex not in $M$, and the suspensions in the before and after states are both filtered.
\end{definition}

Note that when $k = n$, an extended bistellar operation is simply a bistellar operation.

The filtration of $(\dagger)$ then gives a setting to iteratively apply the following:

\begin{lemma} \label{extend} If $B$ is a triangulated ball, then there is a triangulation of $B\times I$ 
in which there is a subcomplex of the form $\Sigma B$, with the vertices $v_+$ and $v_-$
having positive and negative $I$ coordinates, respectively.
\end{lemma}

\noindent{\sc Proof.}  Triangulate $B \times \{1, 0, -1\}$ with three copies of the triangulation on $B$, and $\partial B \times I$ with two copies of one of the usual triangulations of a product of a complex with a 1-simplex induced from an ordering of the vertices of the complex and an ordering of the vertices of the 1-simplex.  The result is a cell complex with exactly two cells which are not simplexes, but which have triangulated boundaries, one consisting of points with positive $I$ coordinate, the other consisting of points with negative $I$ coordinate.  Triangulate each of the untriangulated cells by starring.  Plainly the resulting triangulation has the desired property. $\Box$\medskip

We are now in the position to state and prove our main theorem:

\begin{thm} \label{main}  If $M$ is a filtered manifold, any two triangulations of $M$ are related by a finite sequence of extended bistellar operations.
\end{thm}

\noindent{\sc Proof.}   We proceed by induction on the codimension of $M_k$, the lowest dimensional stratum on which the two triangulations differ.  (Note that $k$ is necessarily positive, since $M_0$ as a finite set of points admits a unique triangulation.)

For the base case of codimension zero, notice that by the extension of Casali's lemma, it follows that for a fixed triangulation of $M_{n-1}$, every triangulation of $M$ which restricts to it can be reached by a sequence of bistellar moves on simplexes in $M\setminus M_{n-1}$ (which are a special case of extended bistellar moves).  

Now, suppose for $i \leq j$ we have shown that all triangulations which agree on $M_{n-i}$ can be obtained from each other by extended bistellar operations on simplexes in strata of codimension less than $i$.  We wish to show that the same holds for $i = j+1$.

Suppose $T$ and $T^\prime$ are triangulations which agree on $M_{n-j-1}$.  If they agree on $M_{n-j}$ we are done by the induction hypothesis.  If they do not, then they necessarily differ on $M_{n-j}\setminus M_{n-j-1}$  If we ignore the ambient filtered manifold $M$ and regard $M_{n-j}$ as a filtered manifold of dimension $n-j$, by the base case of the induction, there is a sequence of bistellar moves in $M_{n-j}\setminus M_{n-j-1}$ which takes $T|_{M_{n-j}}$ to
$T^\prime|_{M_{n-j}}$.

Now, consider the first of these moves.  It is replacing a subcomplex of $M_{n-j}$ of the form
$A\ast \partial B$, for $A$ a simplex of $M_{n-j}\setminus M_{n-j-1}$ and $B$ a simplex not in $M$ with $\partial A \ast B$.  Both of 
$A\ast \partial B$ and $\partial A \ast B$ are triangulated balls (of dimension $n-j$), $j$-fold iterated application of Lemma \ref{extend} then gives a triangulation of the filtered submanifold with boundary $[A\ast \partial B] \times I^{n-j}$ which has $\Sigma^{n-j}[A\ast \partial B]$ as a subcomplex a filtered iterated suspension of $[A\ast \partial B]$ and extends to a triangulation $T^{\prime \prime}$ of $M$, agreeing with $T$ on $M_{n-j}$.  By the inductive hypothesis,
$T^{\prime \prime}$ can be reached from $T$ by a sequence of extended bistellar moves.  The extended bistellar move replacing  $\Sigma^{n-j}[A\ast \partial B]$ with  $\Sigma^{n-j}[\partial A\ast B]$, then performs the first bistellar move in $M_{n-j}\setminus M_{n-j-1}$.

Replacing $T$ with $T^{\prime \prime}$ and iterating the construction just given, eventually reaches a triangulation agreeing with $T^\prime$ on $M_{n-j}$, and one more application of the induction hypothesis completes the proof. $\Box$\medskip

Now, observe that as in the case of Pachner's \cite{Pa87} result, for any given dimension of manifold (here filtered manifold) there are only finitely many combinatorial types of moves needed to relate any two triangulations of the same manifold.  In Pachner's case there were $n-1$ inverse pairs (counting a move of the same combinatorial type as its inverse as an inverse pair) needed for $n$-manifolds.  Here we need $n^2-n$ inverse pairs of moves -- the extensions of the bistellar moves on each stratum.

The key to obtaining this result was the existence of a finite family of local models for neighborhoods of triangulated cells which admitted a standard way of extending any triangulation of the cell which would allow moves to be extended.  We would like to extend the result to more general stratifications.  

We are not able to consider general stratified spaces using only Lemma \ref{BetterC'sLemma} as a tool -- indeed Lemma \ref{BetterC'sLemma} applies to the most commonly considered class of stratified spaces, the CS spaces of Siebenmann \cite{Si72}, only in the very special cases.  We thus define a class of stratified spaces which includes those of interest for the construction of TQFTs with defects:

\begin{definition} A {\em starkly stratified space} is a PL space $X$ equipped with a filtration

\[ X_0 \subset X_1 \subset \ldots \subset X_{n-1} \subset X_n = X \]

\noindent satisfying

\begin{enumerate}
\item There is a triangulation $\cal T$ of $X$ in which each $X_k$ is a subcomplex.
\item For each $k = 1,\ldots n$ $X_k \setminus X_{k-1}$ is a(n open) $k$-manifold.
\item If $C$ is a connected component of $X_k \setminus X_{k-1}$, then $\cal T$ restricted
to $\bar{C}$ gives $\bar{C}$ the structure of a combinatorial manifold with boundary.
\item For each combinatorial ball $B_k$ with $\accentset{\circ}B_k \subset X_k \setminus X_{k-1}$, $\accentset{\circ}B_k$ admits a closed neighborhood $N$ given inductively as a cell complex as follows (although we require $B_k$ to be a combinatorial ball, the triangulation is then ignored):

\noindent $N = N_n$, where $N_m$ for $k \leq m \leq n$ is given inductively by 

\[ N_k = B_k\] 

\noindent and 

\[ N_{\ell+1} = N_\ell \cup \bigcup_{v \in S_\ell} L(v) \ast v \]

\noindent for $S_\ell$ a finite set of points in $X\setminus X_\ell$, andl $L(\cdot)$ a function on $S_\ell$
valued in 

\[ \{ L | L\; \mbox{\rm is a combinatorial ball and}\; B_k \subset L \subset N_\ell  \}\]

We will call such a neighborhood of the interior of a combinatorial ball lying in a stratum of the same dimension a {\em stark neighborhood}.

\end{enumerate}
\end{definition}

Note in the last condition, the points $v$ are only required to lie in $X\setminus X_\ell$, not in
$X_{\ell+1}\setminus X_\ell$.

In the definition, we deliberately ignored the triangulation on $B_k$, but we trivially have the following lemma:

\begin{lemma} \label{starkextension}
If $B_k$ is a combinatorial ball with $\accentset{\circ}B_k \subset X_k \setminus X_{k-1}$, and $N$ is any stark neighborhood of $\accentset{\circ}B_k$, then a triangulation $\cal T$ of $B_k$ extends canonically to a triangulation of $N$, which we denote ${\cal T}\uparrow N$, by iterated coning.
\end{lemma}

It is now evident what the appropriate combinatorial moves are for transforming triangulations of starkly stratified spaces.

\begin{definition} An {\em extended bistellar operation} in a starkly stratified space is the combinatorial operation which replaces $[A \ast \partial B]\uparrow N$ with $[\partial A \ast B]\uparrow N$, where $[A \ast \partial B]$ (resp. $[\partial A \ast B]$ denote the triangulation of a combinatorial $k$-ball with interior in $X_k\setminus X_{k-1}$ as the join of a $j$-simplex $A$ and the boundary of a $k-j$-simplex $B$ (resp. the join of the boundary of the $j$-simplex $A$ and the $k-j$-simplex $B$).
\end{definition}

Note that this agrees with our previous definition of an extended bistellar move in the case of a filtered manifold, since in that case stark neighborhoods and $n-k$-fold iterated filtered suspensions coincide.

Modifying the proof of Theorem \ref{main} by replacing iterated suspensions with the canonical extension to a stark neighborhood now carries a stronger result:

\begin{thm}  If $X$ is a starkly stratified space, any two triangulations of $X$ are related by a finite sequence of extended bistellar operations.
\end{thm}

Plainly this does {\em not} give a finite set of moves that suffice for all starkly stratified spaces of fixed dimension -- this is easily seen by considering the family of two-dimensional starkly stratified spaces given as joins of a 1-simplex and a set of $n$ points for $n = 1, 2, \ldots$ stratified by 
letting $X_0$ be the endpoints of the 1-simplex and $X_1$ be the union of the boundaries of the all the $n$ resulting 2-simplexes.

If, however, we restrict the number of distinct structure-preserving homeomorphism-types of stark neighborhood to be finite, as was the case for filtered manifolds, then for starkly stratified spaces with those specified types of stark neighborhoods as local models, only finitely many distinct types of combinatorial moves are needed.

One interesting case where this restriction can be made are stark stratifications of the form

\[ \emptyset \subset K \subset \Sigma \subset S^3 \] 

\noindent where $K$ is a classical knot and $\Sigma$ is a Seifert surface.  In this case five inverse pairs of moves suffice:  3-dimensional bistellar operations, suspensions of 2-dimensional bistellar operations in $\Sigma$ (with vertices in $S^3\setminus \Sigma$), and the double suspension of 1-dimensional bistellar operations in $K$, in which the first supension has one vertex in $\Sigma$ and the other in $S^3\setminus \Sigma$, and the second has both vertices in $S^3\setminus \Sigma$.

Another case in which the stratifications are naturally stark and in which a finite family of moves suffice are stratified spaces arising from
an immersion $f:N_{n-1} \rightarrow M_n$ which is completely regular in the sense that at all multiple points, the branches meet in general position (cf. \cite{He81}).  $M$ is then naturally stratified by the multiplicity of points with respect to $f$, and it is evident from the local models provided by complete regularity that having triangulated to move from the smooth to the PL setting, every $\accentset{\circ}B_k \subset M_k \setminus M_{k-1}$ admits a stark neighborhood, and moreover the pattern of iterated joins depends only on the codimension.

In future work the authors plan to apply the results herein to the explicit construction of TQFTs from various categorical and algebraic data.   It is also hoped to refine the present results to find families of combinatorial moves which respect flag-like triangulation, that is triangulations of a stratified space in which each the intersection of a simplex with the closure of a stratum is a face of the simplex.  Such triangulations are used in the simplicial formulation of intersection cohomology and have been used to construct TQFTs with defects in codimension-2 by the second author \cite{Ye92}.

\noindent {\sc Department of Mathematics, Kansas State University, Manhattan, Kansas 66506}

\noindent {\tt dyetter@math.ksu.edu}

\noindent {\tt crane@math.ksu.edu}
\end{document}